\documentclass[12pt,twoside,final]{amsart}

\usepackage{amssymb}
\usepackage{times,a4wide}

\theoremstyle{plain}
\newtheorem*{theorem*}{Theorem}

\newtheorem{theorem}{Theorem}[section]

\newtheorem*{proposition*}{Proposition}

\newtheorem*{corollary*}{Corollary}

\newtheorem*{lemma*}{Lemma}

\theoremstyle{definition}

\newtheorem*{remark*}{Remark}

\theoremstyle{definition}

\newtheorem*{definition*}{Definition}
\newcommand{\nc}{\newcommand}

\newcommand{\C}{{\mathbb C}}
\newcommand{\Z}{{\mathbb Z}}

\newcommand{\D}{{\mathbb D}}
\newcommand{\R}{{\mathbb R}}

\newcommand{\T}{{\partial \mathbb D}}

\newcommand{\clos}{\operatorname{clos}}

\newcommand{\eps}{\varepsilon}
\newcommand{\vp}{\varphi}
\nc{\bea}{\begin{eqnarray}}
\nc{\eea}{\end{eqnarray}}
\nc{\beqa}{\begin{eqnarray*}}
\nc{\eeqa}{\end{eqnarray*}}
\nc{\Hi}{H^{\infty}}
\nc{\loi}{\ell^{\infty}}
\nc{\NL}{N^+\vert \Lambda}
\nc{\hf}{{\mathcal H}_{\phi}}
\nc{\liL}{\lambda\in\Lambda}
\nc{\nn}{\nonumber}

\newenvironment{proof*}{\vskip 2mm\noindent {}}{$\blacksquare$ \vskip 2mm}
\numberwithin{equation}{section}

\renewcommand{\Im}{\operatorname{Im}}

\newcommand{\sinc}{\operatorname{sinc}}

\title
{Reverse Carleson Measures in Hardy spaces}

\author{Andreas Hartmann,
Xavier Massaneda, Artur Nicolau, \& Joaquim Ortega-Cerd\`a}

\address{Laboratoire de Math\'ematiques Pures de Bordeaux,
Universit\'e de Bordeaux, 351 cours de la Lib\'eration,
33405 Talence, France}

\address{Departament de Matem\`atica Aplicada i An\`alisi,
Universitat  de Barcelona, Gran Via 585, 08007-Bar\-ce\-lo\-na, Spain}

\address{Departament de Matem\`atiques,
Universitat Aut\`onoma de Barcelona, 08193-Bellaterra, Spain}

\address{Departament de Matem\`atica Aplicada i An\`alisi,
Universitat  de Barcelona, Gran Via 585, 08007-Bar\-ce\-lo\-na, Spain}

\thanks{First author supported by ANR FRAB: ANR-09-BLAN-0058-02 and
CTP network ``Analyse math\'ematique et applications''.
Other authors supported by the Generalitat de Catalunya (grants 2009 SGR 01303
and 00420) and the spanish Ministerio de Ciencia e Innovaci\'on (projects
MTM2011-27932-C02-01 and MTM2011-0374).}

\email{hartmann@math.u-bordeaux.fr,
xavier.massaneda@ub.edu, 
\newline
artur@mat.uab.cat, jortega@ub.edu}

\date{\today}



\begin{document}

\begin{abstract}
We give a necessary and sufficient condition for a measure $\mu$ in the closed
unit disk to be a reverse Carleson measure for Hardy spaces. This extends a
previous result of Lef\`evre, Li, Queff\'elec and Rodr\'{\i}guez-Piazza
\cite{LLQR}. We also provide a simple example showing that the analogue
for the Paley-Wiener space does not hold. This example can be
generalised to model spaces associated to one-component inner functions.
\end{abstract}

\maketitle

\section{Introduction}

For $1\le p<\infty$ let $H^p$ be the Hardy space on the unit disk $\D$ equipped with its usual
norm
\[
 \| f\|_p=\left(\sup_{r < 1} \int_0^{2\pi} |f(re^{i\theta})|^p\;
\frac{d\theta}{2\pi}\right)^{1/p}\ .
\]
Denote by $M_+({\D})$ the set of positive, finite Borel
measures supported on ${\D}$, and let $\mu\in
M_+({\D})$. A well known theorem by Carleson (see
\cite[Chap.I Th. 5.6]{Gar}) states that $H^p$  embeds into
$L^p({\D},\mu)$:
\begin{eqnarray}\label{carleson}
 \|f\|_{L^p(\D,\mu)}\lesssim \|f\|_p,\quad f\in H^p,
\end{eqnarray}
if and only if $\mu$ satisfies the 
Carleson condition: there exists $C>0$ such that for all
arcs $I$ in $\T$
\begin{eqnarray}\label{carlcond}
 \mu(S_I)\leq C |I|,
\end{eqnarray}
where $S_I=\{z\in\overline{\D}:1-|I|\le |z|\le 1, z/|z| \in I\}$ is
the usual Carleson window. This theorem has been extended to
several other spaces, like Bergman, Fock, model spaces etc., and we
refer the reader to the huge bibliography on this topic for further
information.

Note that $H^p$ contains a dense set of continuous functions for which
the embedding \eqref{carleson} still makes sense when the measure has a part
supported on the boundary. Then \eqref{carlcond} implies that the restriction
of the measure $\mu$ to the boundary has to be absolutely continuous
with respect to Lebesgue measure and with bounded Radon-Nikodym
derivative. 
It is thus possible to consider more generally positive, finite
Borel measures supported on the closed unit disk: $M_+(\overline{\D})$.

Here, we are interested in  reverse Carleson inequalities
$\|f\|_p\lesssim \|f\|_{L^p(\overline{\D},\mu)}$, $f \in
C(\overline{\D})\cap H^p(\D)$, $1<p<\infty$. In \cite{LLQR} Lef\`evre et
al.\  proved that when $\mu$ is already a Carleson measure these
hold if and only it there exists $C>0$ such that for all arcs
$I\subset\T$
\[
 \mu(S_I)\geq C |I|.
\]

Our elementary proof actually
shows that the reverse inequalities hold without the Carleson condition.
It turns out that the interesting part of the measure has
to be supported on the boundary, while the part supported in the disk can
be dropped.

The embedding problem is closely related with the \emph{reproducing kernel
thesis}: if the embedding holds on the reproducing kernels, then it actually
holds for every function. We also show that the reproducing kernel thesis holds
for the reverse Carleson embedding.

Finally, we provide a simple example showing that the analogous reproducing
kernel thesis for the reverse embedding in the Paley-Wiener space does not hold.
The construction can be generalised to model spaces associated to one-component
inner functions.

We shall use the following standard notation:   $f \lesssim g$ means that there
is a constant $C$ independent of the relevant variables such that $f\le C g$,
and $f\simeq g$ means that $f\lesssim g$ and $g\lesssim f$.

\section{Main result}

For $1<p<\infty$ and $\lambda\in\D$ consider the reproducing kernel in $H^p$
\[
 k_{\lambda}(z)=\frac{1}{1-\overline{\lambda}z},\quad z\in\D,
\]
and its normalised companion
\[
 K_{\lambda}:=\frac{k_{\lambda}}{\|k_{\lambda}\|_p}\ .
\]
A standard computation shows that $\|k_\lambda\|_p\simeq
(1-|\lambda|)^{-1/p'}$, where $1/p+1/p'=1$.

Our main result reads as follows.

\begin{theorem}
Let $1<p<\infty$ and let $\mu \in M_+(\overline{\D})$.
Then the following assertions are equivalent:
\begin{enumerate}
\item There exists $C_1>0$ such that for every 
function $f\in H^p\cap C( \overline{ \D ) }$,
\[
\displaystyle{\int_{\overline{\D}}|f|^p d\mu}
\geq C_1 \|f\|_p^p\ ,
\]

\item There exists $C_2>0$ such that for every $\lambda\in\D$,
\[
\int_{\overline{\D}}|K_{\lambda}|^pd\mu
\geq C_2\ ,
\]

\item There exists $C_3>0$ such that for every arc $I\subset\T$,
\[
\mu(S_I)\geq C_3 |I|\ .
\]

\item 
There exists $C_4>0$ such that the Radon-Nikodym derivative of
$\mu|_{\T}$ with respect to the length measure is bounded below by $C_4$.
\end{enumerate}
\end{theorem}

Observe that in this theorem we
do not require absolute continuity
of the restriction $\mu|_{\T}$.
Still, if we want to extend (1) to the entire $H^p$-space, then, in order that
$\int_{\overline{\D}}|f|^pd\mu$ makes sense for every function in
$H^p$, we need to impose absolute continuity
on $\mu|_{\T}$.
Note that the integral $\int_{\overline{\D}}|f|^pd\mu$
can be infinite for certain $f\in H^p$ when the Radon-Nikodym derivative
of $\mu|_{\T}$ is not bounded.

\begin{proof}
(1) $\Rightarrow$ (2) is clear.

(3) $\Rightarrow$ (4). Take $h>0$ so that $|I|/h$ is a large integer $N$ and
consider the modified Carleson window
\[
 S_{I,h}=\{z\in \overline{\D} : 1-h\leq |z|\leq 1,\; z/|z| \in I\}\ .
\]
Split $I$ into $N$ subarcs $I_k$ such that $|I_k|=h$ (and hence
$S_{I_k,h}=S_{I_k}$). Then
\[
 \mu(S_{I,h})=\mu(\bigcup_{k=1}^N S_{I_k,h})=
\sum_{k=1}^N \mu(S_{I_k,h})\geq C_3 \sum_{k=1}^N |I_k| = C_3 |I|.
\]
Now, for every open set $O$
in $\overline{\D}$  for which $I\subset O$ there exists $h>0$ such that
$S_{I,h}\subset O$. Since $\mu\in M_+(\D^-)$ is outer regular (see
\cite[Theorem 2.18]{Ru}) we thus have
\[
 \mu(I)=\inf_{I\subset O \text{ open in  } \overline{\D}}\mu(O)\ge\inf_{h>0}\mu(S_{I,h})
 \geq C_3 |I|.
\]
We deduce that the Lebesgue measure on $\T$ denoted by $m$ is absolutely
continuous with respect to the restriction of $\mu$ to $\T$ and that the
corresponding Radon-Nikodym derivative of $\mu$ is bounded below by
$C_3$. In particular one can choose $C_4 = C_3$.

$(4) \Rightarrow (1)$ Clearly, for all $f\in H^p$,
\[
 \int_{\overline{\D}}|f|^p d\mu\ge \int_{\T}|f|^p d\mu \geq C_4 \int_{\T}|f|^pdm
 =C_4 \|f\|_p^p
\]
(in particular, one can choose $C_1=C_4$).

(2) $\Rightarrow$ (3).
By hypothesis, integrating over $S_{I,h}$ with respect to  area measure $dA$ on $\overline{\D}$ we get
\[
 C_2  |I|\times h \leq \int_{S_{I,h}}\int_{\overline{\D}} |K_{\lambda}|^pd\mu dA(\lambda)
 \simeq \int_{\overline{\D}} \int_{S_{I,h}}  \frac{(1-|\lambda|^2)^{p/p'}}{|1-\overline{\lambda}z|^p} dA(\lambda) d\mu(z).
\]

Set
\[
 \varphi_h(z)=\frac{1}{h}\int_{S_{I,h}}
 \frac{(1-|\lambda|^2)^{p/p'}}{|1-\overline{\lambda}z|^p}dA(\lambda)
 =\frac{1}{h}\int_{S_{I,h}}
 \frac{(1-|\lambda|^2)^{p-1}}{|1-\overline{\lambda}z|^p}dA(\lambda),
\]
so that the previous estimate becomes
\begin{equation}\label{ineq}
 \int_{\overline{\D}} \varphi_h(z) d\mu(z)\gtrsim |I|\ .
\end{equation}

We claim that
\[
 \lim_{h\to 0}\varphi_h(z)
 \left\{
 \begin{array}{ll}
 \simeq 1 & \text{if }z\in \overline{I} ,\\
 =0 & \text{otherwise}.
 \end{array}
 \right.
\]
Indeed, if $z\notin \overline{I}$, then there are $\delta,h_0 >0$
such that for every $0<h<h_0$ and for every $\lambda\in S_{I,h}$, we
have $|1-\overline{\lambda}z|\ge\delta>0$, and the result follows
from the estimate
\[
 0\le \vp_h(z)=\frac{1}{h}\int_{S_{I,h}}
 \frac{(1-|\lambda|^2)^{p-1}}{|1-\overline{\lambda}z|^p}dA(\lambda)
 \le \frac{1}{\delta^p}\frac{|I|\times h}{h}\times (2h)^{p-1}\lesssim h^{p-1}.
\]
Suppose now that $z=e^{i\theta_0}\in \overline{I}$.
Let $h\le |I|$, then
setting $\lambda=(1-t) e^{i\theta}$ for $\lambda\in S_{I,h}$
we have
\beqa
 \varphi_h(z)&=&\frac{1}{h}\int_{S_{I,h}}
 \frac{(1-|\lambda|^2)^{p-1}}{|1-\overline{\lambda}z|^p}dA(\lambda)
 \geq\frac{1}{h}\int_{e^{i\theta}\in I}\int_0^h\frac{t^{p-1}}{|e^{i\theta_0}-(1-t)e^{i\theta}|^p}(1-t)dtd\theta\\
 &\gtrsim&\frac{1}{h}\int_0^h\int_{|\theta-\theta_0|\le t, e^{i\theta}\in I}
 \frac{t^{p-1}}{|\theta-\theta_0|^p+t^p}d\theta dt\\
 &\ge& \frac{1}{h}\int_0^h\int_{|\theta-\theta_0|\le t,e^{i \theta}\in I}\frac{t^{p-1}}{2t^p} d\theta  dt.
\eeqa
Since $0\le t\le h\le |I|$ and $z=e^{it}\in\overline{I}$,
the set $\{e^{i\theta}:|\theta-\theta_0|\le t, e^{i\theta}\in I\}$
contains an interval of length at least $t/2$, we get
\beqa
 \varphi_h(z)&\gtrsim&\frac{1}{h}\int_0^h\frac{t}{2}\times\frac{t^{p-1}}{2t^p}dt\simeq 1.
\eeqa

On the other hand, integrating in polar coordinates, we get
\beqa
 \varphi_h(z)&=&\frac{1}{h}\int_{S_{I,h}}
 \frac{(1-|\lambda|^2)^{p-1}}{|1-\overline{\lambda}z|^p}dA(\lambda)
  =\frac{1}{h}\int_{1-h}^1(1-r^2)^{p-1}\int_I\frac{1}{|1-re^{i(\theta-\theta_0)}|^p}
 d\theta rdr\\
 &\lesssim& \frac{1}{h}\int_0^h t^{p-1}\frac{1}{t^{p/p'}} dt
 \simeq 1.
\eeqa Hence $\varphi_h$ converges pointwise to a function comparable
to $\chi_{\overline{I}}$, and $\varphi_h$ is uniformly bounded in
$h$. Now, from \eqref{ineq} and by dominated convergence we finally
deduce that
\[
 \mu(\overline{I})=\int_{\D^-} \chi_{\overline{I}} d\mu \simeq \int_{\overline{\D}}
 \lim_{h\to 0}\varphi_h(z)d\mu(z)=
 \lim_{h\to 0}\int_{\overline{\D}}\varphi_h(z)d\mu(z)
 \gtrsim |I|\ .
\]
\end{proof}

\begin{remark*}
 The following example shows that the reproducing kernel thesis fails
for the reverse Carleson inequalities in the Paley-Wiener space $PW_\pi$,
the space of Fourier transforms of square integrable functions on $[-\pi,\pi]$. In Section 2 we will show how it can be adapted to any model space associated to a one-component inner function.

Consider the sequence $S=\{x_n\}_{n\in\Z\setminus\{0\}}$, where
\[
 x_n=
 \begin{cases}
  n+1/8\ &\textrm{if $n$ is even}\\
  n-1/8\ &\textrm{if $n$ is odd.}
 \end{cases}
\]
By the Kadets-Ingham theorem (see e.g.\ \cite[Theorem D4.1.2]{Nik}) $S$
would be a minimal sampling sequence if we added the point $0$. Since $S$ is not
sampling the
discrete measure $\mu:=\sum_{n\neq 0}\delta_{x_n}$
does not satisfy the reverse inequality $\|f\|_{L^2(\R)}\lesssim \|f\|_{L^2(\mu)}$, $f\in PW_\pi$.

Let us see that, on the other hand, the $\mu$-norm of the normalised
reproducing kernels
\[
 K_{\lambda}(z)=c_\lambda
\sinc(\pi(z-\lambda))=c_\lambda\frac{\sin(\pi(z-\lambda))}{
\pi(z-\lambda)},\qquad c_\lambda^2\simeq (1+ |\Im\lambda|)e^{-2\pi |\Im \lambda|},
\]
is uniformly bounded from below. If $\lambda$ is such that $|\Im
\lambda|>1$ then
$|\sin(\pi(x_n-\lambda))|\simeq e^{\pi|\Im \lambda|}$, and hence
\[
\int_{\C} |K_\lambda(x)|^2 d\mu(x)
 =\sum_{n\neq 0} c_{\lambda}^2\left|\frac{\sin(\pi(x_n-\lambda))}{\pi(x_n-\lambda)}
 \right|^2\simeq
\sum_{n\neq 0} \frac
{|\Im\lambda|}{|x_n-\lambda|^2}\simeq 1.
\]
It is thus enough to consider points $\lambda\in\C$ with $|\Im\,
\lambda|\le 1$. Let $x_{n_0}$ be the point of $S$ closest to
$\lambda$; then there is $\delta>0$, independent of $\lambda$, such
that
\[
 \int_{\C} |K_\lambda(x)|^2 d\mu(x)
 =\sum_{n\neq 0} |K_\lambda (x_n)|^2
 \geq \left|\frac{\sin(\pi(x_{n_0}-\lambda))}{\pi(x_{n_0}-\lambda)}\right|^2\ge \delta\ .
\]
It is interesting to point out that $\mu$ is a Carleson measure for $PW_\pi$, since $S$ is in a strip and separated.
\end{remark*}

\section{Failure in other model spaces}

The previous construction can be generalised to certain model spaces in the disk.
The model space associated to an inner function $\Theta$ is  $K_\Theta=
H^2\ominus \Theta H^2$, and the reproducing kernel corresponding to $\lambda\in\D$ is given by
\[
 k_{\lambda}^\Theta(z)=\frac{1-\overline{\Theta(\lambda)}\Theta(z)}{1-\overline{\lambda}z},
 \quad z\in\D.
\]
A particular class of model spaces is given by the so-called \emph{one-component} inner functions, those for which the sub-level set
$\{z\in\D:|\Theta(z)|<\eps\}$ is connected for some $0<\eps<1$ .

The Paley-Wiener space corresponds, after a conformal mapping of $\D$ into the upper half-plane, to the inner function $\Theta_{2\pi}(z)=e^{i2\pi z}$. More precisely
$K_{\Theta_{2\pi}}=e^{i\pi z} PW_\pi$.

Here we show the following result.

\begin{theorem}\label{CLSsituation}
If $\Theta$ is a one-component inner function, then the reverse reproducing
kernel thesis does not hold in $K_\Theta$.
\end{theorem}

We refer the reader to \cite{BFGHR} for sufficient conditions for reverse
Carleson measures in model spaces.

Let $\sigma(\Theta)$ denote the spectrum of $\Theta$, that is, the set of $\zeta\in\overline{\D}$ such that
 $\liminf\limits_{z\to\zeta, z\in \D} |\Theta(z)|=0$.
For one-component inner functions the set $\T\setminus \sigma(\Theta)$ is a countable union of arcs where $\Theta$ is analytic (and on which the argument of $\Theta$
increases by $2\pi$). Moreover, for
any $|\alpha|=1$,
\[
 E_\alpha : =\{\zeta\in \T\setminus \sigma(\Theta) : \Theta(\zeta)=\alpha\}
\]
is countable and the system $(K^\Theta_{\zeta_n})_{\zeta_n\in E_\alpha}$ is an orthonormal basis of $K_\Theta$, a so-called Clark basis (see \cite{Cl}, and \cite[Section 4]{BaDy}
for the material needed here).
For such $\zeta\in\T\setminus \sigma(\Theta)$ the reproducing kernel is defined as
\[
 k_{\zeta}^\Theta(z)=\frac{1-\overline{\Theta(\zeta)}\Theta(z)}{1-\overline{\zeta}z}
 =\zeta\overline{\Theta(\zeta)}\frac{\Theta(\zeta)-\Theta(z)}{\zeta-z},\quad z\in\D.
\]
Its norm is $\sqrt{|\Theta'(\zeta)|}$, so that the corresponding normalised reproducing kernel is
\[
 K_{\zeta}^\Theta:=\frac{k_{\zeta}^\Theta}{\|k_{\zeta}^\Theta\|_2}=\frac{k_{\zeta}^\Theta}{\sqrt{|\Theta'(\zeta)|}}\ .
\]

With these elements we follow the scheme of the Paley-Wiener case to prove
Theorem \ref{CLSsituation}.

\begin{proof}
Pick the Clark basis $(K^\Theta_{\zeta_n})_{n\ge 0}$ for $\alpha=1$ and
set
\beqa
 \xi_n=
 \left\{\begin{array}{ll}
 \zeta_n & \text{if $n\neq 1$} \\
 \xi_1 & \text{if $n=1$},
 \end{array}
 \right.
\eeqa where we choose $\xi_1 $ sufficiently close  to $\zeta_1$ (and
in particular different from $\zeta_n$, $n\neq 1$) but different
from $\zeta_1$, implying in particular $\langle
K^\Theta_{\xi_1},K^\Theta_{\zeta_n}\rangle\neq 0$ for every $n$, so
that  $(K^\Theta_{\xi_n})_{n\ge 0}$ is an unconditional basis (see
\cite{BaDy}; it is actually not far from being orthogonal). It will
be clear from the proof below how close to $\zeta_1$ we have to
choose $\xi_1$.



We now consider the measure
\[
 \mu:=\sum_{n>0} \|k_{\xi_n}^\Theta\|_2^{-2} \delta_{\xi_n} =\sum_{n>0}|\Theta'(\xi_n)|^{-1}\delta_{\xi_n}
\]
where we have taken away the very first point $\xi_0$, so that
$(K_{\xi_n}^\Theta)_{n > 0}$ is an incomplete family. Notice that this
is a perturbation of the Clark measure $\sigma=\sum_{n\geq 0} \|k_{\zeta_n}^\Theta\|^{-2} \delta_{\zeta_n}$ with
one mass point deleted.
Thus $\mu$ is not a reverse Carleson measure
since
there are functions vanishing in all the points $\xi_n$, $n>0$, but
not in $\xi_0$.\\

Let us check that the reverse reproducing kernel thesis fails, which, in view of the above,
amounts to find a  
$\delta>0$ such that $\|K^\Theta_z\|_{L^2(\mu)}\ge \delta$ for every $z\in \D$.
Note that
\begin{equation}\label{normKz}
 \|K^\Theta_z\|^2_{L^2(\mu)}
 =\sum_{n\ge 1}\frac{1}{|\Theta'(\xi_n)|} |K^\Theta_z(\xi_n)|^2=\sum_{n\ge 1}
 |\langle K^\Theta_z,K^\Theta_{\xi_n}\rangle|^2,
\end{equation}
which are just the generalised Fourier coefficients of $K^\Theta_z$ in $K^\Theta_{\xi_n}$, $n\ge 1$.

Let us introduce the following function
\[
 \varphi(z):=|\langle K^\Theta_{\zeta_0},K^\Theta_z\rangle|^2
 =\left|\frac{\Theta(\zeta_0)-\Theta(z)}{\zeta_0-z}\right|^2 \frac{1}{|\Theta'(\zeta_0)|}
 \frac{1-|z|^2}{1-|\Theta(z)|^2},\quad z\in\D.
\]
By the Cauchy-Schwarz inequality $\vp(z)\le 1$ for all $z\in\D$.
Also,
since $\|K^\Theta_{\zeta_0}\|_2=\|K^\Theta_z\|_2=1$, the only way to get
$\vp(z)=1$ is that $K^\Theta_z=\alpha K^\Theta_{\zeta_0}$, $|\alpha|=1$, i.e.\ $z=\zeta_0$.

Since $\zeta_0$ is not in the spectrum, there is a closed neighbourhood
$C$ of $\zeta_0$ in $\overline{\D}$ on which $\Theta$ is analytic, which implies that
$\varphi$ is continuous on $C$. We suppose $C$ small enough
that it does not contain any other $\zeta_k$, $k\neq 0$, nor $\xi_1$.

Introduce the sets
\[
 U_{\delta}:=\{z\in C :|z-\zeta_0|<\delta\}
\]
and define
\[
 \psi(\delta):=\sup_{z\notin U_{\delta}} \varphi(z)
\]

\emph{Claim:} For sufficiently small $\delta$ the function $\psi(\delta)$
is decreasing,  with $\psi(0)=1$ and $\psi(\delta)<1$ for $\delta>0$.

We postpone the proof of the claim and proceed now to prove that
$\|K^\Theta_z\|_{L^2(\mu)}\gtrsim 1$. Pick $\delta >0$ sufficiently
small such that $\psi(\delta) < 1$. We consider two cases.

Assume first that $z\notin U_{\delta}$. Pick $0 < \varepsilon <
 1- \psi(\delta)$.
Since $\{\zeta_0\}\cup\{\xi_k\}_{k\ge 1}$ gives rise to
a perturbation of the orthonormal Clark basis $(K_{\zeta_n})_{n\ge
0}$, it suffices to choose $\xi_1$ close enough to $\zeta_1$ so that
there is $0<\eta<\epsilon$ such that for every $f\in K_\Theta$ (see
\cite{BaDy})
\[
 (1-\eta) \|f\|_2^2 \le |\langle f,K_{\zeta_0}\rangle|^2 + \sum_{n\ge 1}|\langle f,K_{\xi_n}\rangle|^2
 \le (1+\eta) \|f\|_2^2.
\]
Then, by \eqref{normKz}
\begin{align*}
 \|K_z^\Theta\|_{L^2(\mu)}^2 &=\sum_{n\ge 1}|\langle K^\Theta_z,K^\Theta_{\xi_n}\rangle|^2
 = |\langle K^\Theta_z,K^\Theta_{\zeta_0}\rangle|^2+ \sum_{n\ge 1}|\langle K^\Theta_z,K^\Theta_{\xi_n}\rangle|^2
 -|\langle K^\Theta_z,K^\Theta_{\zeta_0}\rangle|^2
 \\
 &\ge (1-\eta)\|K^\Theta_z\|_2^2-\vp(z)
 \ge (1-\eta)-(1-\epsilon)=\epsilon-\eta>0
\end{align*}

Assume now that $z\in U_{\delta}\subset C$. We will check that on
this set it suffices to consider only two terms of the sum
$\varphi_1(z)  = |\langle K_z^\Theta , K_{\xi_1}^\Theta\rangle|^2$
and $\varphi_2(z) = |\langle K_z^\Theta, K_{\zeta_2}^\Theta
\rangle|^2$.
It is here that we need that $\xi_1$ is a small perturbation of $\zeta_1$ which
is ``not harmonic'' with $\zeta_1$, meaning that $ |\langle K_{\zeta_2}^\Theta,K_{\xi_1}^\Theta\rangle|^2
\neq 0$. Indeed $\varphi_1$ and $\varphi_2$ are continuous functions on
the compact set $\overline U_{\delta}$. Since $U_{\delta}\subset C$, we have
$\varphi_2(z)=0$, $z\in \overline U_{\delta}$, if and only if $z=\zeta_0$. Now
$\varphi_1(\zeta_0)>0$ so that by a continuity argument we conclude that
$\varphi_1(z)+\varphi_2(z)$ is strictly bounded away from 0 for $z\notin U_{\delta}$,
which concludes the proof.

\emph{Proof of the Claim.} It is clear that $\psi(\delta)$
is decreasing and $\psi(0)=1$.

We prove now that $\psi(\delta)<1$ for $\delta>0$.
Indeed, suppose not, then there is a sequence $(z_n)_n\subset \D\setminus
U_\delta$ such that
$\vp(z_n)=|\langle K^\Theta_{\zeta_0},K^\Theta_{z_n}\rangle|^2\to 1$ as $n\to\infty$.
We can also assume that $z_n\to \zeta\in\clos(\T\setminus U_\delta)$.
Now $(K^\Theta_{z_n})_n$ is a bounded family, and by the Alaoglu theorem
it admits a weak convergent subsequence, which in order not
to overcharge notation, we can suppose to be also indexed by $n$.
Let $f$ be a weak limit of this sequence so that
$|\langle K_{\zeta_0}^\Theta,f\rangle|=1$. It is also clear that $\|f\|_2=1$.
From the same observation
as above we can deduce $f=\alpha K^\Theta_{\zeta_0}$, $|\alpha|=1$ (in fact,
every weak convergent subsequence has $K^\Theta_{\zeta_0}$ as weak limit).
In particular, by the weak convergence, for every $f\in K_\Theta$,
\begin{eqnarray}\label{weakcv}
 f(z_n)\sqrt{\frac{1-|z_n|^2}{1-|\Theta(z_n)|^2}}
 = \langle f,K_{z_n}^\Theta\rangle\to\langle f,K_{\zeta_0}^\Theta\rangle=
 \frac{f(\zeta_0)}{\sqrt{|\Theta'(\zeta_0)|}}.
\end{eqnarray}

Observe that $K_\Theta$ contains continuous functions (by a result of
Aleksandrov continuous functions in $K_\Theta$ form actually a dense set in $K_\Theta$,
see \cite[p.186]{CMR}).

Now, if there are two continuous functions $f_1$ and $f_2$ in $K_\Theta$
such that the vectors $(f_1(\zeta),f_1(\zeta_0))$ and
$(f_2(\zeta),f_2(\zeta_0))$
are linearly independent, then we can deduce from \eqref{weakcv}
that necessarily, first
\[
\frac{1-|z_n|^2}{1-|\Theta(z_n)|^2} \to \frac{1}{{|\Theta'(\zeta_0)|}}
\]
and then
\[
 f_1(\zeta)= f_1(\zeta_0)\quad\text{and}\quad  f_2(\zeta)=f_2(\zeta_0)
\]
which is not possible unless $\zeta=\zeta_0$.

Let us prove that if $\zeta\ne \zeta_0$ then there are two such
functions $f_1,f_2$. We start by taking two linearly independent
continuous functions $h_1,h_2\in K_\Theta$. It may happen that
$(h_1(\zeta),h_1(\zeta_0))$ and $(h_2(\zeta),h_2(\zeta_0))$ are
linearly independent and then we are done. If they are linearly
dependent, then we can find a linear combination $f$ of $h_1$ and
$h_2$ which is not identically $0$ and such that
$f(\zeta)=f(\zeta_0)=0$. Consider the backward shift operator
$S^*f(z)=\frac{f(z)-f(0)}z$ and recall that $S^*K_\Theta\subset
K_\Theta$. Observe that if moreover $f(0)=0$ then also
$S^*f(\zeta)=S^*f(\zeta_0)=0$. Hence, after sufficiently many
applications of $S^*$ we can suppose that $f(0)\neq 0$,
$f(\zeta)=f(\zeta_0)=0$, and, renormalising, that $f(0)=1$.

Then $g=S^*f$ is continuous in $K_\theta$ and
takes two different values $g(\zeta)=-\overline{\zeta}$ and
$g(\zeta_0)=-\overline{\zeta_0}$. Set now $h={S^*}^2f$ which takes
the values $h(\zeta)=-\overline{\zeta}^2-\overline{\zeta}h'(0)$ and
$h(\zeta_0)=-\overline{\zeta_0}^2-\overline{\zeta_0}h'(0)$.
Then either the vectors  $(g(\zeta),g(\zeta_0))$ and
$(h(\zeta),h(\zeta_0))$ are linearly independent (and we are done) or
they are not, in which case the solution of the linear dependence gives
$\zeta=\zeta_0$.
\end{proof}

\end{document}